\documentclass[11pt]{article} 
\usepackage{latexsym} 
\newcommand{\eproof}{\mbox{\ }\hfill $\Box$ \par \vskip 10pt}

\newtheorem{Theorem}{Theorem}[section] 
 
\newtheorem{prop}[Theorem]{Proposition}

\begin{document}

\title{Weighted $L^p$ decay estimates of solutions to the
wave equation with a potential}

\author{{\sc  Fernando Cardoso and Georgi Vodev}}

\date{} 
 
\maketitle 

\begin{abstract}
\noindent
We obtain large time decay estimates on weighted $L^p$ spaces, $2<p<+\infty$,
for solutions to the wave equation with real-valued potential
$V(x)=O(\langle x\rangle^{-2-\delta_0})$, $\delta_0>0$. 
\end{abstract}

\setcounter{section}{0}
\section{Introduction and statement of results}

Let $V\in L^\infty({\bf R}^3)$ be a real-valued function satisfying
$$\left|V(x)\right|\le C\langle x\rangle^{-2-\delta_0},
\quad \forall x\in {\bf R}^3,\eqno{(1.1)}$$
with constants $C>0$ and $\delta_0>0$ not necessarily small, 
where $\langle x\rangle=(1+|x|^2)^{1/2}$. 
Denote by $G_0$ and $G$ the self-adjoint realizations of the
operators $-\Delta$ and $-\Delta+V(x)$ on $L^2({\bf R}^3)$.
We suppose that $G$ has no eigenvalues, which in turn implies
$G\ge 0$. Moreover, under the assumption (1.1), $G$
has no strictly positive resonances (e.g. see \cite{kn:GV}, \cite{kn:V}).

Given any $a>0$ denote by $\chi_a\in C^\infty({\bf R})$, 
$\chi_a\ge 0$, a function supported in the interval
$[a,+\infty)$, $\chi_a=1$ on $[a+1/2,+\infty)$.
It is well known that the solutions to the free wave equation
satisfy the following dispersive estimate 
$$\left\|
G_0^{-\alpha}e^{it\sqrt{G_0}}\right\|_{L^{p'}\to
L^p}\le C|t|^{-\alpha},\quad t\neq 0,\eqno{(1.2)}$$
for every  $2\le p<+\infty$, where $1/p+1/p'=1$, $\alpha=1-2/p$.
Hereafter, given $1\le p\le +\infty$, $L^p$ denotes the space
$L^p({\bf R}^3)$. 
It turns out that a better decay is possible to get in weighted
$L^p$ spaces. Namely, we have the following estimate (see the appendix):
$$\left\|\langle x\rangle^{-\sigma\alpha}
G_0^{-\alpha}e^{it\sqrt{G_0}}\chi_a(\sqrt{G_0})
\langle x\rangle^{-\sigma\alpha}\right\|_{L^{p'}\to
L^p}\le C|t|^{-\alpha(1+\sigma)},\quad |t|\ge 1,\eqno{(1.3)}$$
for every $a>0$, $\sigma\ge 0$, 
$2\le p<+\infty$, where $1/p+1/p'=1$, $\alpha=1-2/p$.
 The purpose of this note is to prove
an analogue of (1.3) for the operator $G$. Our main result is the following

\begin{Theorem}
Assume (1.1) fulfilled.
Then, for every $a>0$, $2\le p<+\infty$, $0<\sigma<\delta_0$,
 the following estimate holds
$$\left\|\langle x\rangle^{-\sigma\alpha}
G^{-\alpha}e^{it\sqrt{G}}\chi_a(\sqrt{G})\langle x\rangle^{-\sigma\alpha}
\right\|_{L^{p'}\to
L^p}\le C\left(|t|^{-1-\sigma}\log(1+|t|)\right)^\alpha,
\quad |t|\ge 1,\eqno{(1.4)}$$
 where $1/p+1/p'=1$, $\alpha=1-2/p$, while for $\sigma\ge\delta_0$ we have
$$\left\|\langle x\rangle^{-\sigma\alpha}
G^{-\alpha}e^{it\sqrt{G}}\chi_a(\sqrt{G})\langle x\rangle^{-\sigma\alpha}
\right\|_{L^{p'}\to
L^p}\le C_\epsilon|t|^{-\alpha(1+\delta_0-\epsilon)},
\quad |t|\ge 1,\quad\forall 0<\epsilon\ll 1.\eqno{(1.5)}$$
\end{Theorem}

{\it Remark.} It follows from (1.5) that for potentials 
$V(x)=O_N\left(\langle x
\rangle^{-N}\right)$,$\forall N\gg 1$, we have the estimate
$$\left\|\psi
G^{-\alpha}e^{it\sqrt{G}}\chi_a(\sqrt{G})\psi
\right\|_{L^{p'}\to
L^p}\le C_N|t|^{-N},
\quad |t|\ge 1,\quad\forall N\gg 1,\eqno{(1.6)}$$
for every $2<p<+\infty$ and every function $\psi\in C_0^\infty({\bf R}^3)$.\\

The estimate (1.4) with $\chi_a\equiv 1$, $\sigma=0$, $2\le p\le 4$,
 and without the logarithmic term in the RHS, was
proved in \cite{kn:BS} for potentials
$V(x)=O\left(\langle x\rangle^{-3-\varepsilon_0}\right)$, 
$\varepsilon_0>0$, and later on extended in \cite{kn:GV} to 
non-negative potentials $V(x)=O\left(\langle x
\rangle^{-2-\varepsilon_0}\right)$, $\varepsilon_0>0$. 
Recently, in \cite{kn:CuV} an analogue of (1.4) with $\sigma=0$
was obtained for a larger class of short-range potentials.

To prove Theorem 1.1 we follow some ideas from \cite{kn:CV} and \cite{kn:CuV}.
The proof is based on a carefull study of 
the operator-valued function
$$\langle x\rangle^{-\sigma}(G-\lambda^2\pm i0)^{-1}
\langle x\rangle^{-\sigma}-\langle x\rangle^{-\sigma}
(G_0-\lambda^2\pm i0)^{-1}\langle x\rangle^{-\sigma}
:L^1\to L^\infty,
\quad\lambda\ge\lambda_0>0,$$
together with its derivatives (see Proposition 3.2).
This in turn requires sharp estimates for the resolvent of the
perturbed operator as well as of its derivatives on weighted
$L^2$ spaces (see Proposition 2.2).

{\bf Acknowledgements.} A part of this work was carried out while
the first author was visiting the University of Nantes in June, 2004,
under the support of the agreement Brazil-France in Mathematics - Proc.
69.0014/01-5.

\section{Uniform resolvent estimates}

Given any $\lambda\ge 0$, $0<\varepsilon\le 1$, define the free resolvent by
$$R_0(\lambda\pm i\varepsilon)=\left
(G_0-(\lambda\pm i\varepsilon\right)^2)^{-1}:
L^2\to L^2,$$
with kernel
$$[R_0(\lambda\pm i\varepsilon)](x,y)=
\frac{e^{(\pm i\lambda-\varepsilon)|x-y|}}{4\pi|x-y|}.$$
Then the kernel of 
$ R_0^{(k)}=d^kR_0/d\lambda^k$, $k\ge 1$, is given by
 $$[R_0^{(k)}(\lambda\pm i\varepsilon)](x,y)=
\frac{(\pm i)^k}{4\pi}|x-y|^{k-1}e^{(\pm i\lambda-\varepsilon)|x-y|}.$$

\begin{prop} Let $s>-1/2$, $s_1>1/2$, $0\le \sigma \le 1$, $\lambda_0>0$,
 and let $k\ge 1$ be an integer. Then the following estimates hold:
$$\left\|\langle x\rangle^{-s_1}R_0(\lambda\pm i\varepsilon)
\langle x\rangle^{-s_1}\right\|_{L^2\to L^2}\le
 C\lambda^{-1},\quad\lambda\ge \lambda_0,\eqno{(2.1)}$$
$$\left\|\langle x\rangle^{-s}R_0(\lambda\pm i\varepsilon)
\langle x\rangle^{-s_1}\right\|_{L^2\to L^2}\le
 C\lambda^{-1}\varepsilon^{-\max\{1/2-s+\epsilon,0\}},
\quad\lambda\ge\lambda_0,
\eqno{(2.2)}$$
$$\left\|\langle x\rangle^{-k-s} R_0^{(k)}(\lambda\pm i\varepsilon)
\langle x\rangle^{-k-s_1}\right\|_{L^2\to L^2}\le
 C\lambda^{-1}\varepsilon^{-\max\{1/2-s+\epsilon,0\}},
\quad\lambda\ge\lambda_0,\eqno{(2.3)}$$
$$\left\|R_0(\lambda\pm i\varepsilon)
\langle x\rangle^{-s}\right\|_{L^2\to L^\infty}\le
  C\varepsilon^{-\max\{1/2-s+\epsilon,0\}},\quad\lambda\ge 0,
\eqno{(2.4)}$$
$$\left\| R_0^{(1)}(\lambda\pm i\varepsilon)
\langle x\rangle^{-1-s}\right\|_{L^2\to L^\infty}\le
C\varepsilon^{-\max\{1/2-s+\epsilon,0\}},\quad\lambda\ge 0,
\eqno{(2.5)}$$
$$\left\|\langle x\rangle^{-k+1-\sigma} R_0^{(k+1)}(\lambda\pm i\varepsilon)
\langle x\rangle^{-k-1-s}\right\|_{L^2\to L^\infty}\le
C\varepsilon^{-1+\min\{\sigma,s+1/2-\epsilon\}},
\quad\lambda\ge 0,\eqno{(2.6)}$$
$$\left\|\langle x\rangle^{-s}R_0(\lambda\pm i\varepsilon)
\right\|_{L^1\to L^2}\le
  C\varepsilon^{-\max\{1/2-s+\epsilon,0\}},\quad\lambda\ge 0,
\eqno{(2.7)}$$
$$\left\|\langle x\rangle^{-1-s}R^{(1)}_0(\lambda\pm i\varepsilon)
\right\|_{L^1\to L^2}\le
C\varepsilon^{-\max\{1/2-s+\epsilon,0\}},\quad\lambda\ge 0,
\eqno{(2.8)}$$
$$\left\|\langle x\rangle^{-k-1-s} R_0^{(k+1)}(\lambda\pm i\varepsilon)
\langle x\rangle^{-k+1-\sigma}
\right\|_{L^1\to L^2}\le
C\varepsilon^{-1+\min\{\sigma,s+1/2-\epsilon\}},
\quad\lambda\ge 0,\eqno{(2.9)}$$
$\forall\,0<\epsilon \ll 1$, with a constant
$C=C(\epsilon)>0$ independent of $\lambda$ and $\varepsilon$, where the 
estimates (2.4) and (2.7) hold for $s\ge 0$ only. Moreover, the constant
$C$ in (2.2) and (2.3) may depend on $\lambda_0$. 
\end{prop}

{\it Proof.} The estimate (2.1) is well known and that is why we omit 
its proof.  
The estimates (2.2), (2.4), (2.5), (2.7) and (2.8) 
are proved in \cite{kn:CuV} for $s\ge 0$. Here we will provide the proof
of these estimates (exept for (2.4) and (2.7))
 for $-1/2<s<0$. The proof of (2.2) 
in this case is a little bit 
more involved, while the proof of (2.5) and (2.8)
goes in precisely the same way and we 
will present it just for the sake of completeness. 

To prove (2.2) for this range of values of $s$ we will take advantage
of the formula
$$\langle x\rangle^{-s}R_0(\lambda\pm i\varepsilon)
\langle x\rangle^{-s_1}=R_0(\lambda\pm i\varepsilon)
\langle x\rangle^{-s-s_1}+R_0(\lambda\pm i\varepsilon)[-\Delta,
\langle x\rangle^{-s}]R_0(\lambda\pm i\varepsilon)
\langle x\rangle^{-s_1}.$$
Taking into account that $[-\Delta,
\langle x\rangle^{-s}]=\sum_{j=1}^3 O\left(\langle x\rangle^{-s-1}
\right)\partial_{x_j}+ O\left(\langle x\rangle^{-s-2}
\right)$, we obtain from the above represantation
$$\left\|\langle x\rangle^{-s}R_0(\lambda\pm i\varepsilon)
\langle x\rangle^{-s_1}\right\|_{L^2\to L^2}\le\left\|
R_0(\lambda\pm i\varepsilon)
\langle x\rangle^{-s-s_1}\right\|_{L^2\to L^2}$$ 
$$+C\left\|R_0(\lambda\pm i\varepsilon)\langle x\rangle^{-1/2-\epsilon}
\right\|_{L^2\to L^2}\left\|
\langle x\rangle^{-s-1/2+\epsilon}\nabla R_0(\lambda\pm i\varepsilon)
\langle x\rangle^{-s_1}\right\|_{L^2\to L^2}$$
 $$+C\left\|R_0(\lambda\pm i\varepsilon)\langle x\rangle^{-1/2-\epsilon}
\right\|_{L^2\to L^2}\left\|
\langle x\rangle^{-s-3/2+\epsilon} R_0(\lambda\pm i\varepsilon)
\langle x\rangle^{-s_1}\right\|_{L^2\to L^2},$$
for $0<\epsilon\ll 1$. On the other, it is not hard to see that when 
$s\ge 0$ (2.2) implies
$$\left\|\langle x\rangle^{-s}\nabla R_0(\lambda\pm i\varepsilon)
\langle x\rangle^{-s_1}\right\|_{L^2\to L^2}\le
 C\varepsilon^{-\max\{1/2-s+\epsilon,0\}},
\quad\lambda\ge\lambda_0,$$
with a constant $C>0$ depending on $\lambda_0$. 
Now it is easy to see that (2.2) with $-1/2<s<0$ follows from the above
estimates and (2.2) with $s\ge 0$. 

Let us now see that (2.3) follows from (2.2) by induction in $k$. Set
$$-\widetilde\Delta:=-r\Delta r^{-1}=-\partial^2_r+r^{-2}\Delta_{S^2},$$
where $\Delta_{S^2}$ denotes the (positive) Laplace-Beltrami operator
on $S^2:=\{x\in{\bf R}^3:|x|=1\}$, and 
denote by $\widetilde G_0$ the self-adjoint realization of the operator
$-\widetilde\Delta$ on the Hilbert space $H=L^2({\bf R}^+\times S^2,
drdw)$. Clearly, $\widetilde G_0$ is unitary equivalent to $G_0$, so it
suffices to prove (2.3) with $G_0$ and $L^2$ replaced by $\widetilde G_0$
and $H$, respectively. Using the identity
$$-2\widetilde\Delta+
[r\partial_r,-\widetilde\Delta]=0,\eqno{(2.10)}$$
we obtain the following representation for the first derivative of 
the resolvent
$$\widetilde R_0(\lambda\pm i\varepsilon)=(\widetilde G_0-(
\lambda\pm i\varepsilon)^2)^{-1}:H\to H,$$
namely
$$(\lambda\pm i\varepsilon)
\langle r\rangle^{-k-s}\widetilde R_0^{(1)}(\lambda\pm i\varepsilon)
\langle r\rangle^{-k-s_1}=-\langle r\rangle^{-k-s}
\widetilde R_0(\lambda\pm i\varepsilon)
\langle r\rangle^{-k-s_1}$$ $$
+\langle r\rangle^{-k-s}\widetilde R_0(\lambda\pm i\varepsilon)
\partial_rr\langle r\rangle^{-k-s_1}
-\langle r\rangle^{-k-s}r\partial_r
\widetilde R_0(\lambda\pm i\varepsilon)
\langle r\rangle^{-k-s_1}.$$
 Differentiating this identity $k-1$ times with respect to $\lambda$ 
leads to 
$$(\lambda\pm i\varepsilon)
\langle r\rangle^{-k-s}\widetilde R_0^{(k)}(\lambda\pm i\varepsilon)
\langle r\rangle^{-k-s_1}=-2\langle r\rangle^{-k-s}
\widetilde R_0^{(k-1)}(\lambda\pm i\varepsilon)
\langle r\rangle^{-k-s_1}$$ $$
+\langle r\rangle^{-k-s}\widetilde R_0^{(k-1)}(\lambda\pm i\varepsilon)
\partial_rr\langle r\rangle^{-k-s_1}
-\langle r\rangle^{-k-s}r\partial_r
\widetilde R_0^{(k-1)}(\lambda\pm i\varepsilon)
\langle r\rangle^{-k-s_1}.$$ On the other hand, it is easy to see 
(for example, this follows from the estimate (2.22) below obtained in
a more general situation) that (2.3) with $k-1$ implies
$$\left\|\langle r\rangle^{-k-s}r\partial_r
\widetilde R_0^{(k-1)}(\lambda\pm i\varepsilon)
\langle r\rangle^{-k-s_1}\right\|_{H\to H}\le 
C\varepsilon^{-\max\{1/2-s+\epsilon,0\}},
\quad\lambda\ge\lambda_0,$$ 
with a constant $C>0$ depending on $\lambda_0$. 
Therefore, the estimate (2.3) with $k-1$ implies
(2.3) with $k$.

To prove (2.5) for $-1/2<s<0$ observe that we have
$$\left\| R_0^{(1)}(\lambda\pm i\varepsilon)
\langle x\rangle^{-1-s}\right\|_{L^2\to L^\infty}^2\le
\sup_{x\in{\bf R}^3} A_s(x,\varepsilon),$$
where
$$A_s(x,\varepsilon)=
\int_{{\bf R}^3}e^{-2\varepsilon|x-y|}\langle y
\rangle^{-2s-2}dy=\int_{|y|\le |x|/2}+
\int_{|y|\ge |x|/2}$$
 $$\le e^{-\varepsilon|x|}\int_{|y|\le |x|/2}\langle y
\rangle^{-2s-2}dy+C\int_{|y|\ge |x|/2}
e^{-2\varepsilon|x-y|}\langle x-y\rangle^{-2s-2}dy$$
 $$\le Ce^{-\varepsilon|x|}\int_0^{|x|/2}(\rho+1)^{-2s}d\rho+
C\int_{0}^\infty e^{-2\varepsilon\rho}(\rho+1)^{-2s}d\rho$$
 $$\le C'e^{-\varepsilon|x|}(|x|+1)^{-2s+1}+
C'\int_{0}^\infty e^{-2\varepsilon\rho}(\rho^{-2s}+1)d\rho
\le C''\varepsilon^{-1+2s}.$$
The estimate (2.8) is obtained in precisely the same way.

In what follows we will prove (2.6) and (2.9). We have
$$\left\|\langle x\rangle^{-k+1-\sigma} R_0^{(k+1)}(\lambda\pm i\varepsilon)
\langle x\rangle^{-k-1-s}\right\|_{L^2\to L^\infty}^2\le
\sup_{x\in{\bf R}^3} B_{s,\sigma,k}(x,\varepsilon),$$
where
$$ B_{s,\sigma,k}(x,\varepsilon)=\langle x\rangle^{-2\sigma-2k+2}
\int_{{\bf R}^3}|x-y|^{2k}e^{-2\varepsilon|x-y|}\langle y
\rangle^{-2s-2k-2}dy$$ $$=\langle x\rangle^{-2\sigma-2k+2}\int_{|y|\le |x|/2}+
\langle x\rangle^{-2\sigma-2k+2}\int_{|y|\ge |x|/2}$$
 $$\le C\langle x\rangle^{-2\sigma-2k+2}
|x|^{2k}e^{-\varepsilon|x|}\int_{|y|\le |x|/2}\langle y
\rangle^{-2s-2k-2}dy+C\int_{{\bf R}^3}
e^{-2\varepsilon|z|}\langle z\rangle^{-2s-2}dz$$
 $$\le C\varepsilon^{2\sigma-2}+
C_\epsilon\varepsilon^{-2\max\{1/2-s+\epsilon,0\}},$$
 $\forall 0<\epsilon\ll 1$, uniformly in $x$. The estimate (2.9) 
is obtained in precisely the same way.
\eproof

Define the perturbed resolvent by
$$R(\lambda\pm i\varepsilon)=\left
(G-(\lambda\pm i\varepsilon\right)^2)^{-1}:L^2\to L^2,$$
and denote 
$R^{(k)}(\lambda\pm i\varepsilon):=d^kR(\lambda\pm i\varepsilon)/d\lambda^k$, 
$k\ge 1$. Let $k_0\ge 0$ be the bigest integer strictly less than
$\delta_0$, and set $\delta'_0=\delta_0-k_0>0$.

\begin{prop} Under the assumption (1.1),
 for every $\lambda_0>0$, $s>-1/2$, $s_1>1/2$, $\lambda\ge\lambda_0$,
$0<\varepsilon\le 1$, the following estimates hold:
$$\left\|\langle x\rangle^{-s}R(\lambda\pm i\varepsilon)
\langle x\rangle^{-s_1}\right\|_{L^2\to L^2}\le
 C_\epsilon\lambda^{-1}\varepsilon^{-\max\{1/2-s+\epsilon,0\}}, 
\eqno{(2.11)}$$
$$\left\|\langle x\rangle^{-k-s} R^{(k)}(\lambda\pm i\varepsilon)
\langle x\rangle^{-k-s_1}\right\|_{L^2\to L^2}\le
 C_\epsilon\lambda^{-1}\varepsilon^{-\max\{1/2-s+\epsilon,0\}},\quad
k=1,...,k_0+1, \eqno{(2.12)}$$
$$\left\|\langle x\rangle^{-k_0-2-s}R^{(k_0+2)}(\lambda\pm i\varepsilon)
\langle x\rangle^{-k_0-2-s}\right\|_{L^2\to L^2}\le
 C_\epsilon\lambda^{-1}\varepsilon^{-1+\min\{s+1/2-\epsilon,\delta'_0\}}, 
\eqno{(2.13)}$$
$\forall\,0<\epsilon\ll 1$, 
where the constant $ C_\epsilon>0$ may depend also on $\lambda_0$.
\end{prop}

{\it Proof.} Clearly, it suffices to prove (2.11) and (2.12) for $1/2<s_1\le
(1+\delta'_0)/2$. We are going to take advantage of the identity
$$\langle x\rangle^{-j-s}R(\lambda\pm i\varepsilon)
\langle x\rangle^{-j-s_1}\left(1+K_j(\lambda\pm i\varepsilon)\right)=
\langle x\rangle^{-j-s}R_0(\lambda\pm i\varepsilon)
\langle x\rangle^{-j-s_1},\eqno{(2.14)}$$
where $0\le j\le k_0+1$, and the operator
$$K_j(\lambda\pm i\varepsilon)=\langle x\rangle^{j+s_1}V
R_0(\lambda\pm i\varepsilon)\langle x\rangle^{-j-s_1}$$
takes values in the compact operators in ${\cal L}(L^2)$. By (2.1),
we have
$$\|K_j(\lambda\pm i\varepsilon)\|_{L^2\to L^2}\le C\lambda^{-1}\le 1/2,
\eqno{(2.15)}$$
for $\lambda\ge\lambda_0$, $0\le\varepsilon\le 1$, with some
$\lambda_0>0$. Hence, for these values of $\lambda$ we have
$$\|\left(1+K_j(\lambda\pm i\varepsilon)\right)^{-1}
\|_{L^2\to L^2}\le C,\eqno{(2.16)}$$
with a constant
$C>0$ independent of $\lambda$ and $\varepsilon$. Moreover, since the operator
$G$ has no strictly positive real resonances, it is easy to see that in fact
(2.16) holds for every $\lambda_0>0$ with a constant $C>0$ depending on
$\lambda_0$. Then (2.11) follows from (2.2), (2.14) and (2.16).

Differentiating (2.14) $j$ times, we get
$$\langle x\rangle^{-j-s} R^{(j)}(\lambda\pm i\varepsilon)
\langle x\rangle^{-j-s_1}\left(1+K_j(\lambda\pm i\varepsilon)\right)
=\langle x\rangle^{-j-s} R_0^{(j)}(\lambda\pm i\varepsilon)
\langle x\rangle^{-j-s_1}$$ 
 $$-\sum_{\nu=0}^{j-1}\beta_{\nu,j}
\langle x\rangle^{-j-s} R^{(\nu)}
(\lambda\pm i\varepsilon)\langle x\rangle^{-\nu-s_1}V
\langle x\rangle^{\nu+s_1} R_0^{(j-\nu)}(\lambda\pm i\varepsilon)
\langle x\rangle^{-j-s_1}.\eqno{(2.17)}$$
Now it is easy to see that (2.12) follows by induction 
from (2.16) and (2.17) combined with (2.3) and (2.11). 
 
To prove (2.13) we will proceed in a way similar to that one
 in \cite{kn:CuV}. 
Denote by $\widetilde G$ the self-adjoint realization of the operator
$-\widetilde\Delta+V$ on the Hilbert space $H$ introduced above. 
Clearly, $\widetilde G$ is unitary equivalent to $G$, so it
suffices to prove (2.13) with $G$ and $L^2$ replaced by $\widetilde G$
and $H$, respectively. Using (2.10) 
we obtain the following representation for the first derivative of 
the resolvent
$$\widetilde R(\lambda\pm i\varepsilon)=(\widetilde G-(
\lambda\pm i\varepsilon)^2)^{-1}:H\to H,$$
namely
$$(\lambda\pm i\varepsilon)
\langle r\rangle^{-2-k_0-s}\widetilde R^{(1)}(\lambda\pm i\varepsilon)
\langle r\rangle^{-2-k_0-s}=-\langle r\rangle^{-2-k_0-s}
\widetilde R(\lambda\pm i\varepsilon)
\langle r\rangle^{-2-k_0-s}$$ $$
+\langle r\rangle^{-2-k_0-s}\widetilde R(\lambda\pm i\varepsilon)
\partial_rr\langle r\rangle^{-2-k_0-s}
-\langle r\rangle^{-2-k_0-s}r\partial_r
\widetilde R(\lambda\pm i\varepsilon)
\langle r\rangle^{-2-k_0-s}$$
 $$+\langle r\rangle^{-2-k_0-s}\widetilde R(\lambda\pm i\varepsilon)
\partial_rrV\widetilde R(\lambda\pm i\varepsilon)
\langle r\rangle^{-2-k_0-s}$$ 
$$-\langle r\rangle^{-2-k_0-s}\widetilde R(\lambda\pm i\varepsilon)
Vr\partial_r\widetilde R(\lambda\pm i\varepsilon)
\langle r\rangle^{-2-k_0-s}$$  
 $$+\langle r\rangle^{-2-k_0-s}\widetilde R(\lambda\pm i\varepsilon)
V\widetilde R(\lambda\pm i\varepsilon)
\langle r\rangle^{-2-k_0-s}.\eqno{(2.18)}$$
Differentiating (2.18) $k_0+1$ times with respect to $\lambda$ leads to 
$$(\lambda\pm i\varepsilon)
\langle r\rangle^{-2-k_0-s}\widetilde R^{(k_0+2)}(\lambda\pm i\varepsilon)
\langle r\rangle^{-2-k_0-s}=-2\langle r\rangle^{-2-k_0-s}
\widetilde R^{(k_0+1)}(\lambda\pm i\varepsilon)
\langle r\rangle^{-2-k_0-s}$$ $$
+\langle r\rangle^{-2-k_0-s}\widetilde R^{(k_0+1)}(\lambda\pm i\varepsilon)
\partial_rr\langle r\rangle^{-2-k_0-s}
-\langle r\rangle^{-2-k_0-s}r\partial_r
\widetilde R^{(k_0+1)}(\lambda\pm i\varepsilon)
\langle r\rangle^{-2-k_0-s}$$
$$+\sum_{\nu=0}^{k_0+1}\alpha_\nu
\langle r\rangle^{-2-k_0-s}\widetilde R^{(k_0+1-\nu)}(\lambda\pm i\varepsilon)
\partial_rrV\widetilde R^{(\nu)}(\lambda\pm i\varepsilon)
\langle r\rangle^{-2-k_0-s}$$ 
$$-\sum_{\nu=0}^{k_0+1}\alpha_\nu
\langle r\rangle^{-2-k_0-s}\widetilde R^{(k_0+1-\nu)}(\lambda\pm i\varepsilon)
Vr\partial_r\widetilde R^{(\nu)}(\lambda\pm i\varepsilon)
\langle r\rangle^{-2-k_0-s}$$ 
 $$+\sum_{\nu=0}^{k_0+1}\alpha_\nu\langle r\rangle^{-2-k_0-s}
\widetilde R^{(k_0+1-\nu)}(\lambda\pm i\varepsilon)
V\widetilde R^{(\nu)}(\lambda\pm i\varepsilon)
\langle r\rangle^{-2-k_0-s}.\eqno{(2.19)}$$
By (2.19), we obtain
$$|\lambda\pm i\varepsilon|\left\|\langle r\rangle^{-2-k_0-s}
\widetilde R^{(k_0+2)}(\lambda\pm i\varepsilon)
\langle r\rangle^{-2-k_0-s}\right\|\le
 C\left\|\langle r\rangle^{-2-k_0-s}
 \widetilde R^{(k_0+1)}(\lambda\pm i\varepsilon)
\langle r\rangle^{-2-k_0-s}\right\|$$
 $$+C\left\|b_{s+k_0+1}(r)\partial_r
\widetilde R^{(k_0+1)}(\lambda\pm i\varepsilon)
\langle r\rangle^{-2-k_0-s}\right\|
+C\left\|\langle r\rangle^{-2-k_0-s}
\widetilde R^{(k_0+1)}(\lambda\pm i\varepsilon)\partial_r
b_{s+k_0+1}(r)\right\|$$ 
 $$+C\sum_{\nu=0}^{k_0+1}\left\|\langle r\rangle^{-2-k_0-s}
\widetilde R^{(k_0+1-\nu)}(\lambda\pm i\varepsilon)
\partial_rb_{s_0+k_0+1-\nu}(r)\right\|
\left\|\langle r\rangle^{-s_0-\nu}
\widetilde R^{(\nu)}(\lambda\pm i\varepsilon)
\langle r\rangle^{-2-k_0-s}\right\|$$
 $$+C\sum_{\nu=0}^{k_0+1}\left\|\langle r\rangle^{-2-k_0-s}
\widetilde R^{(\nu)}(\lambda\pm i\varepsilon)
\langle r\rangle^{-s_0-\nu}\right\|
\left\|b_{s_0+k_0+1-\nu}(r)\partial_r
\widetilde R^{(k_0+1-\nu)}(\lambda\pm i\varepsilon)
\langle r\rangle^{-2-k_0-s}\right\|$$
 $$+C\sum_{\nu=0}^{k_0+1}\left\|\langle r\rangle^{-2-k_0-s}
\widetilde R^{(\nu)}(\lambda\pm i\varepsilon)
\langle r\rangle^{-s_0-\nu}\right\|
\left\|\langle r\rangle^{-s_0-k_0-1+\nu}
\widetilde R^{(k_0+1-\nu)}(\lambda\pm i\varepsilon)
\langle r\rangle^{-2-k_0-s}\right\|,\eqno{(2.20)}$$
where $b_s(r)=r\langle r\rangle^{-1-s}$, $s_0=\delta'_0/2$, and
$\|\cdot\|$ denotes the norm on ${\cal L}(H)$.

Given any $f\in H$, the function $u=\widetilde R
(\lambda\pm i\varepsilon)f$ satisfies the equation
$$\left(-\partial^2_r+r^{-2}\Delta_{S^2}+V-(\lambda\pm i\varepsilon)^2
\right)u=f,$$
so we have
$$\left(-\partial^2_r+r^{-2}\Delta_{S^2}+V-(\lambda\pm i\varepsilon)^2
\right)b_s(r)u=b_s(r)f+
[-\partial_r^2,b_s(r)]u.$$
Integrating by parts yields, $\forall\gamma>0$,
$$\left\|\partial_r(b_s(r)u)\right\|^2_H\le 
\left(C_1+|\lambda\pm i\varepsilon|^2\right)
\left\|\langle r\rangle^{-s}u\right\|^2_H$$
 $$+\left|\left\langle b_s(r)f+
[-\partial_r^2,b_s(r)]u,b_s(r)u\right\rangle_H\right|$$
 $$\le 
\left(C_2+|\lambda\pm i\varepsilon|^2\right)
\left\|\langle r\rangle^{-s}u\right\|^2_H+C_3\|b_s(r)f\|^2_H
 +\gamma^2\left\|r\langle r\rangle^{-1}
[-\partial_r^2,b_s(r)]u\right\|^2_H$$
 $$\le 
\left(C_4+|\lambda\pm i\varepsilon|^2\right)
\left\|\langle r\rangle^{-s}u\right\|^2_H+C_3\|b_s(r)f\|^2_H
 +O(\gamma^2)\left\|b_s(r)\partial_ru\right\|^2_H.$$
Hence,
$$\left\|b_s(r)\partial_ru\right\|_H\le
\left\|\partial_r\left(b_s(r)u\right)\right\|_H
+C\left\|\langle r\rangle^{-s}u\right\|_H$$
 $$\le 
\left(C+|\lambda\pm i\varepsilon|\right)
\left\|\langle r\rangle^{-s}u\right\|_H+C\|b_s(r)f\|_H
 +O(\gamma)\left\|b_s(r)\partial_ru\right\|_H,$$
which, after taking $\gamma$ small enough, gives
$$\left\|b_s(r)\partial_ru\right\|_H\le
C\left(1+|\lambda\pm i\varepsilon|\right)
\left\|\langle r\rangle^{-s}u\right\|_H+C\|b_s(r)f\|_H.
\eqno{(2.21)}$$
By (2.21) we get, for $j=0,1,...$,
$$\left\|b_{s+j}(r)\partial_r
\widetilde R^{(j)}(\lambda\pm i\varepsilon)
\langle r\rangle^{-2-k_0-s}\right\|\le C(1+|\lambda\pm i\varepsilon|)
\left\|\langle r\rangle^{-j-s}
\widetilde R^{(j)}(\lambda\pm i\varepsilon)
\langle r\rangle^{-2-k_0-s}\right\|$$ $$+C
|\lambda\pm i\varepsilon|
\left\|\langle r\rangle^{-j-s}
\widetilde R^{(j-1)}(\lambda\pm i\varepsilon)
\langle r\rangle^{-2-k_0-s}\right\|,\eqno{(2.22)}$$
where $(\lambda\pm i\varepsilon)\widetilde R^{(-1)}(\lambda\pm i\varepsilon)
:=Id$. 
By (2.20) and (2.22) combined with (2.11) and (2.12) we obtain
$$|\lambda\pm i\varepsilon|\left\|\langle x\rangle^{-2-k_0-s}
R^{(k_0+2)}(\lambda\pm i\varepsilon)
\langle x\rangle^{-2-k_0-s}\right\|_{L^2\to L^2}\le
 C_\epsilon\varepsilon^{-1+\min\{1/2+s-\epsilon,1\}}
+C_\epsilon\varepsilon^{-1+\delta'_0-\epsilon},$$
which clearly implies (2.13).
\eproof

\section{Time decay estimates}

Given parameters $A\gg a>0$, choose a function
$\psi_{a,A}\in C_0^\infty([a,A])$ such that
$$\left|\partial^k_y\psi_{a,A}(y)
\right|\le C_k,\quad\forall k\ge 0,$$
with a constant $C_k>0$ independent of $A$. For $z\in{\bf C}$, set
$$\Phi_{A,z}(t)=\langle x\rangle^{-\sigma z}
G^{-z}e^{it\sqrt{G}}\psi_{a,A}(\sqrt{G})\langle x\rangle^{-\sigma z}-
\langle x\rangle^{-\sigma z}
G_0^{-z}e^{it\sqrt{G_0}}\psi_{a,A}(\sqrt{G_0})
\langle x\rangle^{-\sigma z}.$$

\begin{Theorem} Under the assumption (1.1), for every 
$a>0$, $2\le p\le+\infty$, $A\gg a$, $0<\sigma<\delta_0$, we have
$$\left\|\Phi_{A,\alpha}(t)\right\|_{L^{p'}\to
L^p}\le C\left(|t|^{-1-\sigma}\log A\right)^\alpha,
\quad |t|\ge 1,\eqno{(3.1)}$$
 with a constant
$C>0$ independent of $t$ and $A$, where $1/p+1/p'=1$, $\alpha=1-2/p$.
\end{Theorem}

{\it Proof.} We will first prove (3.1) for $p=+\infty$, $p'=1$. We have
$$\Phi_{A,z}(t)=\int_0^\infty
 e^{it\lambda}\lambda^{1-2z}\psi_{a,A}(\lambda)T(\lambda;\sigma z)d\lambda,
\eqno{(3.2)}$$
where $z\in {\bf C}$, ${\rm Re}\,z=1$, and
$$T(\lambda;\sigma z)=(\pi i)^{-1}\left(T^+(\lambda;\sigma z)-
T^-(\lambda;\sigma z)\right),$$
$$T^\pm(\lambda;\sigma z)=\langle x\rangle^{-\sigma z}\left(
R(\lambda\pm i0)-R_0(\lambda\pm i0)\right)\langle x\rangle^{-\sigma z}
$$ $$=\langle x\rangle^{-\sigma z}
R_0(\lambda\pm i0)VR_0(\lambda\pm i0)\langle x\rangle^{-\sigma z}
+\langle x\rangle^{-\sigma z}R_0(\lambda\pm i0)V
R(\lambda\pm i0)VR_0(\lambda\pm i0)
\langle x\rangle^{-\sigma z}.$$
Denote by $j_0\ge 0$ the bigest integer strictly less than $\sigma$,
and denote $\sigma'=\sigma-j_0>0$.

\begin{prop} Under the assumption (1.1), if $0<\sigma<\delta_0$  
the operator-valued functions $T^\pm(\lambda;\sigma z):L^1\to
 L^\infty$ satisfy the estimates
$$\|\partial_\lambda^jT^\pm(\lambda;
\sigma z)\|_{L^1\to L^\infty}\le C,\quad\lambda\ge\lambda_0,\quad
j=0,1,..., j_0+1,\eqno{(3.3)}$$
$$\|\partial_\lambda^{j_0+1}T^\pm(\lambda_2;\sigma z)-
\partial_\lambda^{j_0+1}T^\pm(\lambda_1;\sigma z)\|_{L^1\to L^\infty}
\le C|\lambda_2-\lambda_1|^{\sigma'},\quad\lambda_2,\lambda_1\ge\lambda_0,
\eqno{(3.4)}$$
$\forall\,\lambda_0>0$ with a constant $C>0$ which may depend on
$\lambda_0$ but is independent of $\lambda$, $\lambda_1$, $\lambda_2$ and
$z$. If $\sigma=\delta_0$ we have (3.3) and (3.4) with $j_0=k_0$
and $\forall\,0<\sigma'<\delta'_0$.
\end{prop}

{\it Proof.} Let first $0<\sigma<\delta_0$. This implies $j_0\le k_0$. 
For every integer $j\ge 0$ we have
$$\partial_\lambda^jT^\pm(\lambda;\sigma z)$$ $$=
\sum_{\nu_1+\nu_2+\nu_3=j}\alpha_{\nu_1,\nu_2,\nu_3}
\langle x\rangle^{-\sigma z}\partial_\lambda^{\nu_1}R_0(\lambda\pm i0)V
\partial_\lambda^{\nu_2}R(\lambda\pm i0)V\partial_\lambda^{
\nu_3}R_0(\lambda\pm i0)\langle x\rangle^{-\sigma z}$$
 $$+\sum_{\mu_1+\mu_2=j}\beta_{\mu_1,\mu_2}
\langle x\rangle^{-\sigma z}\partial_\lambda^{\mu_1}
R_0(\lambda\pm i0)V\partial_\lambda^{\mu_2}
R_0(\lambda\pm i0)\langle x\rangle^{-\sigma z}$$
 $$:=\sum_{\nu_1+\nu_2+\nu_3=j}{\cal A}_{\nu_1,\nu_2,\nu_3}
(\lambda;\varepsilon)
+\sum_{\mu_1+\mu_2=j}{\cal B}_{\mu_1,\mu_2}(\lambda;\varepsilon).
\eqno{(3.5)}$$
Let $j\le j_0+1$. Applying (2.4), (2.5), (2.7), (2.8) with $s>1/2$,
(2.6) with $k=\nu_1-1$ (when $\nu_1\ge 2$), (2.7) with $k=\nu_3-1$ 
(when $\nu_3\ge 2$), and (2.11), (2.12) with $k=\nu_2$, $s>1/2$, we get
$$\left\|{\cal A}_{\nu_1,\nu_2,\nu_3}
(\lambda;\varepsilon)\right\|_{L^1\to L^\infty}\le C,\eqno{(3.6)}$$
with a constant $C>0$ independent of $\lambda$ and $\varepsilon$.
Clearly, a similar estimate holds for 
${\cal B}_{\mu_1,\mu_2}(\lambda;\varepsilon)$, and hence (3.3) follows.

To prove (3.4) it suffices to show that
$$\|\partial_\lambda^{j_0+2}T^\pm(\lambda\pm i\varepsilon;
\sigma z)\|_{L^1\to L^\infty}
\le C\varepsilon^{-1+\sigma'}.\eqno{(3.7)}$$
Indeed, (3.7) implies, $\forall 0<\varepsilon\le 1$,
$$\|\partial_\lambda^{j_0+1}T^\pm(\lambda_2\pm i\varepsilon;\sigma z)-
\partial_\lambda^{j_0+1}T^\pm(\lambda_1\pm i\varepsilon;\sigma z)
\|_{L^1\to L^\infty}
\le C|\lambda_2-\lambda_1|\varepsilon^{-1+\sigma'},$$
$$\|\partial_\lambda^{j_0+1}T^\pm(\lambda_k\pm i\varepsilon;\sigma z)-
\partial_\lambda^{j_0+1}T^\pm(\lambda_k;\sigma z)
\|_{L^1\to L^\infty}\le C\varepsilon^{\sigma'},\quad k=1,2,$$
which yield, $\forall 0<\varepsilon\le 1$,
$$\|\partial_\lambda^{j_0+1}T^\pm(\lambda_2;\sigma z)-
\partial_\lambda^{j_0+1}T^\pm(\lambda_1;\sigma z)\|_{L^1\to L^\infty}
\le C\varepsilon^{\sigma'}\left(2+|\lambda_2-\lambda_1|
\varepsilon^{-1}\right).\eqno{(3.8)}$$
Thus, (3.4) follows from (3.8) by taking $\varepsilon=|\lambda_2-\lambda_1|$.

To prove (3.7) we will make use of (3.5) with $j=j_0+2$. If $j_0<k_0$ we have
$j\le k_0+1$ and hence this case can be treated in the same way as above
to get (3.7) with $\sigma'=1$. Let now $j_0=k_0$. Then, since $\sigma<
\delta_0$, we have $\sigma'<\delta'_0$. Using (2.4), (2.7) with 
$s=1/2+\epsilon$, $0<\epsilon\ll 1$, and (2.13) with 
$s=\delta'_0-1/2-\epsilon$, we obtain
$$\left\|{\cal A}_{0,j_0+2,0}
(\lambda;\varepsilon)\right\|_{L^1\to L^\infty}\le C
\left\|R_0(\lambda\pm i\varepsilon)\langle x\rangle^{-1/2-\epsilon}
\right\|_{L^2\to L^\infty}$$
 $$\left\|\langle x\rangle^{-2-k_0-\delta'_0+
1/2+\epsilon}R^{(k_0+2)}
(\lambda\pm i\varepsilon)\langle x\rangle^{-2-k_0-\delta'_0+
1/2+\epsilon}\right\|_{L^2\to L^2}
\left\|\langle x\rangle^{-1/2-\epsilon}R_0(\lambda\pm i\varepsilon)
\right\|_{L^1\to L^2}$$
 $$\le C_\epsilon\varepsilon^{-1+\delta'_0-2\epsilon}\le
C\varepsilon^{-1+\sigma'},\eqno{(3.9)}$$
provided $\epsilon>0$ is taken small enough. Using (2.5), (2.7) with
$s=1/2+\epsilon$, and (2.12) with $k=k_0+1$, $s_1=1/2+\delta'_0-\epsilon$,
$s=\delta'_0-1/2-\epsilon$, we obtain
 $$\left\|{\cal A}_{1,j_0+1,0}
(\lambda;\varepsilon)\right\|_{L^1\to L^\infty}\le C
\left\|R_0^{(1)}(\lambda\pm i\varepsilon)\langle x\rangle^{-3/2-\epsilon}
\right\|_{L^2\to L^\infty}$$
 $$\left\|\langle x\rangle^{-1-k_0-\delta'_0+
1/2+\epsilon}R^{(k_0+1)}
(\lambda\pm i\varepsilon)\langle x\rangle^{-2-k_0-\delta'_0+
1/2+\epsilon}\right\|_{L^2\to L^2}
\left\|\langle x\rangle^{-1/2-\epsilon}R_0(\lambda\pm i\varepsilon)
\right\|_{L^1\to L^2}$$
 $$\le C_\epsilon\varepsilon^{-1+\delta'_0-2\epsilon}\le
C\varepsilon^{-1+\sigma'},\eqno{(3.10)}$$
and similarly for ${\cal A}_{0,j_0+1,1}(\lambda;\varepsilon)$.
Let now $\nu_1=j_0+2$. Using (2.6) with $k=j_0+1$, $\sigma=\sigma'$,
$s=\delta'_0-1/2-\epsilon$, (2.7) with $s=1/2+\epsilon$, and (2.11)
with $s=s_1=1/2+\epsilon$, we obtain
$$\left\|{\cal A}_{j_0+2,0,0}
(\lambda;\varepsilon)\right\|_{L^1\to L^\infty}\le C
\left\|\langle x\rangle^{-j_0-\sigma'}
R_0^{(j_0+2)}(\lambda\pm i\varepsilon)\langle x\rangle^{-j_0-2-\delta'_0
+1/2+\epsilon}
\right\|_{L^2\to L^\infty}$$
 $$\left\|\langle x\rangle^{-
1/2-\epsilon}R
(\lambda\pm i\varepsilon)\langle x\rangle^{-
1/2-\epsilon}\right\|_{L^2\to L^2}
\left\|\langle x\rangle^{-1/2-\epsilon}R_0(\lambda\pm i\varepsilon)
\right\|_{L^1\to L^2}$$
 $$\le 
C\varepsilon^{-1+\sigma'},\eqno{(3.11)}$$
and similarly for ${\cal A}_{0,0,j_0+2}(\lambda;\varepsilon)$.
Let now $\nu_2\le j_0$, $\nu_1\le j_0+1$, $\nu_3\le
j_0+1$. This implies $\nu_1+\nu_3\ge 2$, $\nu_2+\nu_3\ge 1$, 
$\nu_2+\nu_1\ge 1$. As above we have
$$\left\|{\cal A}_{\nu_1,\nu_2,\nu_3}
(\lambda;\varepsilon)\right\|_{L^1\to L^\infty}\le C
\left\|\langle x\rangle^{-\nu_1+1-\sigma'}
R_0^{(\nu_1)}(\lambda\pm i\varepsilon)\langle x\rangle^{-\nu_1-
\nu_3-\delta'_0+1/2+\epsilon}
\right\|_{L^2\to L^\infty}$$
 $$\left\|\langle x\rangle^{-\nu_2-
1/2-\epsilon}R^{(\nu_2)}
(\lambda\pm i\varepsilon)\langle x\rangle^{-\nu_2-
1/2-\epsilon}\right\|_{L^2\to L^2}$$ $$
\left\|\langle x\rangle^{-\nu_1-\nu_3-\delta'_0+1/2
+\epsilon}R_0^{(\nu_3)}
(\lambda\pm i\varepsilon)\langle x\rangle^{-\nu_3+1-\sigma'}
\right\|_{L^1\to L^2}$$
 $$\le C_\epsilon\varepsilon^{-1+\delta'_0-2\epsilon}\le 
C\varepsilon^{-1+\sigma'}.\eqno{(3.12)}$$
It follows from (3.9)-(3.12) that the first sum in the RHS of (3.5) is
$O\left(\varepsilon^{-1+\sigma'}\right)$. In the same way it is 
easy to see that the second sum satisfies the same bound, and hence
(3.7) follows. When $\sigma=\delta_0$, as above one can show that 
(3.7) holds with any $0<\sigma'<\delta'_0.$
\eproof

We will first consider the case of $0<\sigma<\delta_0$. 
Let $\phi\in C_0^\infty([1/3,1/2])$ be a real-valued function,
$\phi\ge 0$, such that $\int\phi(y)dy=1$. Then, the function
$$T_\epsilon^\pm(\lambda;\sigma z)=\epsilon^{-1}
\int T^\pm(\lambda-y;\sigma z)\phi(y/
\epsilon)dy,\quad 0<\epsilon\ll 1,$$
is smooth with values in ${\cal L}(L^1,L^\infty)$ and, in view of 
(3.3) and (3.4),
satisfies the estimates
$$\|\partial_\lambda^jT_\epsilon^\pm(\lambda;
\sigma z)\|_{L^1\to L^\infty}\le C,\quad
j=0,1,..., j_0+1,\eqno{(3.13)}$$
$$\|\partial_\lambda^{j_0+1}T_\epsilon^\pm(\lambda;\sigma z)-
\partial_\lambda^{j_0+1}T^\pm(\lambda;\sigma z)\|_{L^1\to L^\infty}$$ $$
\le \epsilon^{-1}\int\|\partial_\lambda^{j_0+1}T^\pm(\lambda;\sigma z)-
\partial_\lambda^{j_0+1}T^\pm(\lambda-y;\sigma z)
\|_{L^1\to L^\infty}\phi(y/\epsilon)dy$$
 $$\le C\epsilon^{-1}\int y^{\sigma'}\phi(y/\epsilon)dy=
O(\epsilon^{\sigma'}).\eqno{(3.14)}$$
Let us see that we also have
$$\|\partial_\lambda^{j_0+2}T_\epsilon^\pm(\lambda;\sigma z)
\|_{L^1\to L^\infty}\le  C\epsilon^{-1+\sigma'}.\eqno{(3.15)}$$
Given any $0\le\varepsilon\le 1$, define
$$T_\epsilon^\pm(\lambda\pm i\varepsilon;\sigma z)=\epsilon^{-1}
\int T^\pm(\lambda\pm i\varepsilon-y;\sigma z)\phi(y/\epsilon)dy.$$
In view of (3.7), we have
$$\|\partial_\lambda^{j_0+2}
 T_\epsilon^\pm(\lambda\pm i\varepsilon;\sigma z)\|_{L^1\to L^\infty}\le 
C\varepsilon^{-1+\sigma'},\eqno{(3.16)}$$
with a constant $C>0$ independent of $\lambda$, $\varepsilon$ and
$\epsilon$. On the other hand,
$$\|\partial_\lambda^{j_0+2}T_\epsilon^\pm(\lambda\pm i\varepsilon;\sigma z)-
\partial_\lambda^{j_0+2}
T_\epsilon^\pm(\lambda;\sigma z)\|_{L^1\to L^\infty}$$ $$\le 
\epsilon^{-2}\int\|\partial_\lambda^{j_0+1}
T_\epsilon^\pm(\lambda-y\pm i\varepsilon;\sigma z)-
 \partial_\lambda^{j_0+1}T_\epsilon^\pm(\lambda-y;\sigma z)
\|_{L^1\to L^\infty}|\phi'(y/\epsilon)|
dy$$ $$\le C\varepsilon^{\sigma'}\epsilon^{-2}\int
|\phi'(y/\epsilon)|dy\le 
C\varepsilon^{\sigma'}\epsilon^{-1}.\eqno{(3.17)}$$
By (3.16) and (3.17),
$$\|\partial_\lambda^{j_0+2}
T_\epsilon^\pm(\lambda;\sigma z)\|_{L^1\to L^\infty}\le C
\varepsilon^{\sigma'}(\varepsilon^{-1}+\epsilon^{-1}),$$
which implies (3.15) if we take $\varepsilon=\epsilon$.

Integrating by parts we obtain
$$(it)^{j_0+1}\Phi_{A,z}(t)=\int_0^\infty e^{it\lambda}
\frac{d^{j_0+1}}{d\lambda^{j_0+1}}\left(\lambda^{1-2z}\psi_{a,A}(\lambda)
T(\lambda;\sigma z)\right)d\lambda$$
 $$=\sum_{\nu=0}^{j_0+1}\gamma_\nu
\int_0^\infty e^{it\lambda}
\frac{d^{j_0+1-\nu}}{d\lambda^{j_0+1-\nu}}
\left(\lambda^{1-2z}\psi_{a,A}(\lambda)
\right)\partial_\lambda^\nu T(\lambda;\sigma z)d\lambda$$
 $$=(it)^{-1}\sum_{\nu=0}^{j_0}\gamma_\nu
\int_0^\infty e^{it\lambda}\frac{d}{d\lambda}\left(
\frac{d^{j_0+1-\nu}}{d\lambda^{j_0+1-\nu}}
\left(\lambda^{1-2z}\psi_{a,A}(\lambda)
\right)\partial_\lambda^\nu T(\lambda;\sigma z)\right)d\lambda$$
 $$+\int_0^\infty e^{it\lambda}\lambda^{1-2z}\psi_{a,A}(\lambda)
\partial_\lambda^{j_0+1}T(\lambda;\sigma z)d\lambda:=I_1(t)+I_2(t).$$
In view of (3.3) we have
$$\left\|I_1(t)\right\|_{L^1\to L^\infty}
\le C|t|^{-1}\langle z\rangle^{j_0+2}\log A.\eqno{(3.18)}$$
On the other hand, by (3.14) we get
$$\left\|\int_0^\infty
 e^{it\lambda}\lambda^{1-2z}\psi_{a,A}(\lambda)\left(
\partial_\lambda^{j_0+1}T^\pm(\lambda;\sigma z)-\partial_\lambda^{j_0+1}
T_\epsilon^\pm(\lambda;\sigma z)\right)d\lambda
\right\|_{L^1\to L^\infty}$$ $$\le 
C\epsilon^{\sigma'}\int\lambda^{-1}|\psi_{a,A}(\lambda)|d\lambda
\le C\epsilon^{\sigma'}\log A.\eqno{(3.19)}$$
By (3.15) we get
$$\left\|\int_0^\infty
 e^{it\lambda}\lambda^{1-2z}\psi_{a,A}(\lambda)\partial_\lambda^{j_0+1}
T_\epsilon^\pm(\lambda;\sigma z)d\lambda\right\|_{L^1\to L^\infty}$$
 $$=\left\|t^{-1}\int_0^\infty
 e^{it\lambda}\frac{d}{d\lambda}\left(\lambda^{1-2z}\psi_{a,A}(\lambda)
\partial_\lambda^{j_0+1}
T_\epsilon^\pm(\lambda;\sigma z)\right)d\lambda\right\|_{L^1\to L^\infty}
\le C|t|^{-1}\epsilon^{-1+\sigma'}\langle z\rangle\log A.\eqno{(3.20)}$$
Taking $\epsilon=|t|^{-1}$ we deduce from (3.19) and (3.20),
$$\left\|I_{2}(t)\right\|_{L^1\to L^\infty}
\le C|t|^{-\sigma'}\langle z\rangle\log A.\eqno{(3.21)}$$
By (3.18) and (3.21),
$$\left\|\Phi_{A,z}(t)\right\|_{L^1\to L^\infty}
\le C|t|^{-1-\sigma}\langle z\rangle^{j_0+2}\log A,\eqno{(3.22)}$$
$\forall z\in {\bf C}$, ${\rm Re}\,z=1$. 
On the other hand, we have the trivial estimate
$$\left\|\Phi_{A,z}(t)\right\|_{L^2\to L^2}
\le C,\eqno{(3.23)}$$
$\forall z\in {\bf C}$, ${\rm Re}\,z=0$. Now (3.1) follows from (3.22) and
(3.23) by analytic interpolation.
\eproof

Let $\varphi\in C_0^\infty([1/2,2])$, independent of the parameter $A$.
For $z\in{\bf C}$, set
$$F_{A,z}(t)=\langle x\rangle^{-\sigma z}
G^{-z}e^{it\sqrt{G}}\varphi(\sqrt{G}/A)\langle x\rangle^{-\sigma z}-
\langle x\rangle^{-\sigma z}G_0^{-z}e^{it\sqrt{G_0}}\varphi(\sqrt{G_0}/A)
\langle x\rangle^{-\sigma z}.$$
The following proposition is proved in \cite{kn:CuV} 
for potentials $V(x)=O\left(\langle x\rangle^{-1-\varepsilon_0}\right)$, 
 $\varepsilon_0>0$. 

\begin{prop} For every $2\le p\le +\infty$, $|t|\ge 1$, $A\gg 1$, we have
$$\left\|F_{A,\alpha}(t)
\right\|_{L^{p'}\to L^p}\le C
|t|^{2/p}A^{-2/p},\eqno{(3.24)}$$
with a constant $C>0$ independent of $t$ and
$A$, where $1/p+1/p'=1$, $\alpha=1-2/p$.
\end{prop}

We are going to show now that Theorem 3.1 together with
Proposition 3.3 imply Theorem 1.1.
Choose a function $\varphi\in C_0^\infty([1/2,1])$ such that
$\int\varphi(y)dy=1$, and denote $\varphi_1(y)=y
\varphi(y)$. For $A\gg 1$ we can write
$$\chi_a(y)=\psi_{A,a}(y)+\eta_A(y),\eqno{(3.25)}$$
where
$$\psi_{A,a}(y)=\chi_a(y)A^{-1}\int_{y}^{+\infty}
\varphi(\tau/A)d\tau,$$
$$\eta_A(y)=A^{-1}\int_{0}^{y}\varphi(\tau/A)d\tau
=\int_0^1\varphi_1(sy/A)s^{-1}ds.$$
In view of (3.24), since $p<+\infty$, we have
$$\left\|\langle x\rangle^{-\sigma \alpha}
G^{-\alpha}e^{it\sqrt{G}}\eta_{A}(\sqrt{G})
\langle x\rangle^{-\sigma \alpha}-
\langle x\rangle^{-\sigma \alpha}
G_0^{-\alpha}e^{it\sqrt{G_0}}\eta_{A}(\sqrt{G_0})
\langle x\rangle^{-\sigma \alpha}\right\|_{L^{p'}\to
L^p}$$
 $$=\left\|\int_0^1\left(\langle x\rangle^{-\sigma \alpha}
G^{-\alpha}e^{it\sqrt{G}}\varphi_1(s\sqrt{G}/A)
\langle x\rangle^{-\sigma \alpha}\right.\right.$$ 
 $$\left.\left.
-\langle x\rangle^{-\sigma \alpha}
G_0^{-\alpha}e^{it\sqrt{G_0}}\varphi_1(s\sqrt{G_0}/A)
\langle x\rangle^{-\sigma \alpha}\right)s^{-1}ds
\right\|_{L^{p'}\to L^p}$$
 $$\le C|t|^{2/p}A^{-2/p}\int_0^1s^{2/p -1}ds\le C'|t|^{2/p}A^{-2/p}.
\eqno{(3.26)}$$
Combining (3.25), (3.26) and (3.1) we get
$$\left\|\langle x\rangle^{-\sigma \alpha}
G^{-\alpha}e^{it\sqrt{G}}\chi_a(\sqrt{G})
\langle x\rangle^{-\sigma \alpha}-\langle x\rangle^{-\sigma \alpha}
G_0^{-\alpha}e^{it\sqrt{G_0}}\chi_a(\sqrt{G_0})
\langle x\rangle^{-\sigma \alpha}\right\|_{L^{p'}\to
L^p}$$ $$\le C\left(|t|^{-1-\sigma}\log A\right)^\alpha
+C'|t|^{2/p}A^{-2/p},\eqno{(3.27)}$$
for every $A\gg 1$. Now (1.4) follows from (3.27) by taking
$A=(|t|+1)^k$ with $k>0$ big enough, together with (1.3).
Note finally that (1.5) follows in the same way
 by observing that when $\sigma=\delta_0$
(3.21) holds for every $0<\sigma'<\delta'_0$.  
\eproof

\noindent
{\bf Appendix}\\

In what follows we will derive (1.3) from (1.2).
Denote by $\eta(x,t)$ the characteristic function of the set
$\{|x|\le|t|/4\}$. Indeed, it is easy to see that (1.3) follows from
(1.2), the fact that the operator $\chi_a(\sqrt{G_0})$ is bounded on
$L^p$, $2\le p<+\infty$, and the following estimate
$$\left\|\eta 
G_0^{-\alpha}e^{it\sqrt{G_0}}\chi_a(\sqrt{G_0})\eta\right\|_{L^{p'}\to
L^p}\le C_N|t|^{-\alpha N},\quad |t|\ge 1,\eqno{(A.1)}$$
for every integer $N\ge 1$ and every 
$2\le p\le +\infty$, where $1/p+1/p'=1$, $\alpha=1-2/p$. To prove (A.1)
observe that the kernel, $K_z(x,y;t)$, of the operator
$\eta G_0^{-z}e^{it\sqrt{G_0}}\chi_a(\sqrt{G_0})\eta$, $z\in
{\bf C},{\rm Re}\,z=1$, is given by the oscilatory integral
$$K_z(x,y;t)=\eta(x,t)\eta(y,t)(2\pi)^{-3}\int_{{\bf R}^3}e^{it|\xi|
-i\langle x-y,\xi\rangle}|\xi|^{-2z}\chi_a(|\xi|)d\xi$$
 $$=\eta(x,t)\eta(y,t)(2\pi)^{-3}\int_{{\bf S}^2}
\int_0^\infty e^{i\rho(t
-\langle x-y,w\rangle)}\rho^{2-2z}\chi_a(\rho)d\rho dw$$
 $$=c_N\int_{{\bf S}^2}\eta(x,t)\eta(y,t)(t-\langle x-y,w\rangle)^{-N}
\int_0^\infty e^{i\rho(t
-\langle x-y,w\rangle)}\partial_\rho^N\left(
\rho^{2-2z}\chi_a(\rho)\right)d\rho dw.$$
Hence,
$$|K_z(x,y;t)|\le C'_N|t|^{-N}\int_0^\infty\left|\partial_\rho^N\left(
\rho^{2-2z}\chi_a(\rho)\right)\right|d\rho 
 \le C_N\langle z\rangle^N|t|^{-N},\eqno{(A.2)}$$
for every integer $N\ge 2$. By (A.2) we obtain
$$\left\|\eta 
G_0^{-z}e^{it\sqrt{G_0}}\chi_a(\sqrt{G_0})\eta\right\|_{L^{1}\to
L^\infty}\le C_N\langle z\rangle^N
|t|^{- N},\quad |t|\ge 1,\eqno{(A.3)}$$
for every integer $N\ge 2$, $z\in
{\bf C},{\rm Re}\,z=1$. Now (A.1) follows from (A.3) 
by analytic interpolation.

F. Cardoso, Universidade Federal de Pernambuco,
Departamento de Matem\`atica,
CEP. 50540-740 Recife-Pe, Brazil

e-mail: fernando@dmat.ufpe.br

G. Vodev, Universit\'e de Nantes,
 D\'epartement de Math\'ematiques, UMR 6629 du CNRS,
 2, rue de la Houssini\`ere, BP 92208, 44332 Nantes Cedex 03, France

e-mail: vodev@math.univ-nantes.fr

\end{document}